\documentclass[12pt]{amsart}
\usepackage{amscd,amsmath,amssymb,amsfonts}
\theoremstyle{plain}
\newtheorem{thm}{Theorem}
\newtheorem{lem}[thm]{Lemma}
\newtheorem{cor}[thm]{Corollary}
\newtheorem{prop}[thm]{Proposition}

\theoremstyle{definition}
\newtheorem{defn}[thm]{Definition}
\newtheorem{rmk}[thm]{Remark}

\newtheorem{ex}[thm]{Example}

\numberwithin{thm}{section}
\numberwithin{equation}{section}

\newcommand{\eq}[2]{\begin{equation}\label{#1}#2 \end{equation}}
\newcommand{\ml}[2]{\begin{multline}\label{#1}#2 \end{multline}}
\newcommand{\ga}[2]{\begin{gather}\label{#1}#2 \end{gather}}

\newcommand{\surj}{\twoheadrightarrow}
\newcommand{\inj}{\hookrightarrow}

\newcommand{\codim}{{\rm codim}}

\newcommand{\Pic}{{\rm Pic}}
\newcommand{\Div}{{\rm Div}}
\newcommand{\Hom}{{\rm Hom}}

\newcommand{\Spec}{{\rm Spec \,}}

\newcommand{\sF}{{\mathcal F}}

\newcommand{\sH}{{\mathcal H}}

\newcommand{\sL}{{\mathcal L}}

\newcommand{\sO}{{\mathcal O}}

\newcommand{\sZ}{{\mathcal Z}}
\newcommand{\A}{{\mathbb A}}

\newcommand{\C}{{\mathbb C}}

\newcommand{\G}{{\mathbb G}}

\renewcommand{\P}{{\mathbb P}}
\newcommand{\Q}{{\mathbb Q}}

\newcommand{\Z}{{\mathbb Z}}

\begin{document}

\title[K\"unneth Projectors]{
K\"unneth Projectors for Open Varieties}
\author{Spencer Bloch}
\address{Dept. of Mathematics,
University of Chicago,
Chicago, IL 60637,
USA}
\email{bloch@math.uchicago.edu}

\author{H\'el\`ene Esnault}
\address{Mathematik,
Universit\"at Duisburg-Essen, FB6, Mathematik, 45117 Essen, Germany}
\email{esnault@uni-essen.de}

\date{Feb. 20, 2005}
\dedicatory{To Jacob Murre}
\begin{abstract}
We consider correspondences on smooth quasiprojective varieties $U$. An algebraic cycle inducing the K\"unneth projector 
onto $H^1(U)$ is constructed. Assuming normal crossings at infinity, the existence of relative
motivic cohomology is shown to imply the independence of $\ell$ for traces of open correspondences. 
\end{abstract}
\maketitle
\begin{quote}

\end{quote}

\section{Introduction}

Let $X$ be a smooth, projective algebraic variety over an algebraically closed field $k$, and let
$H^*(X)$ denote a Weil cohomology theory. The existence of algebraic cycles on $X\times X$
inducing as correspondences the various K\"unneth projectors $\pi^i : H^*(X) \to H^i(X)$ is one
of the standard conjectures of Grothendieck, \cite{Kl1}, \cite{Kl2}. It is known in general only for the
cases $i=0,1, 2d-1, 2d$ where $d=\dim X$. The purpose of this note is to consider correspondences
on smooth quasi-projective varieties $U$. In the first section we prove the existence of an
``algebraic'' K\"unneth projector $\pi^1: H^*(U) \to H^1(U)$ assuming that $U$ admits a smooth,
projective completion $X$. The word algebraic is placed in quotes here because in fact the
algebraic cycle on
$X\times U$ inducing $\pi^1$ is not, as one might imagine, trivialized on $(X-U)\times U$. It is
only partially trivialized. This partial trivialization is sufficient to define a class in
$H^{2d-1}_c(U) \otimes H^1(U)$ giving the desired projection. Of course, our cycle on $X\times U$
will be trivialized on $(X\setminus U)\times V$ for $V \subset U$ suitably small nonempty open, but our
method does not in any obvious way yield a full trivialization on $(X\setminus U)\times U$. We finish this
first section with some comments on $\pi^i$ for $i>1$ and some speculation, mostly coming from discussions with A. Beilinson,  on how these ideas
might be applied to study the Milnor conjecture that the Galois cohomology ring of the function
field $H^*(k(X), \Z/n\Z)$ is generated by $H^1$. 

In the last section, we use the existence of relative motivic cohomology \cite{Le} to 
 prove an integrality and independence of $\ell$ result for the trace of an
algebraic correspondence
$\Gamma$ on $U\times U$. We are endebted to G. Laumon for pointing out that one may endeavor to prove
this using results already in the literature(\cite{De}, \cite{P}, \cite{Z}, and \cite{Fu}) by
reduction mod $p$ and composition with a high power of Frobenius. Our objective in what follows is to
show how techniques in motivic cohomology can apply to such questions, at least when the divisor at
infinity has normal crossings. 

When the Zariski closure of the correspondence stabilizes the various
strata $D_I$ at infinity (e.g. when the correspondence is the graph of frobenius) then the trace on
$H^*(U)$ is realized as an alternating sum of traces on $H^*(D_I)$. When in addition all the
intersections with the diagonals are transverse, the contribution to the alternating sum coming from
points lying off $U$ cancels, and the trace on $H^*(U)$ is just the sum of the fixed points on $U$. 

We would like to acknowledge helpful correspondence with A. Beilinson, M. Levine, J. Murre, T. Saito, 
and V. Srinivas. We thank G. Laumon and L. Lafforgue for explaining to us \cite{Fu}.

\section{The first K\"unneth component}\label{sect2}
Let $k$ be an algebraically closed field. We work in the category of algebraic
varieties over $k$. $H^*(X)$ will denote \'etale cohomology with
$\Q_\ell$-coefficients for some $\ell$ prime to the characteristic of $k$. If
$k = \C$, we take Betti cohomology with $\Q$-coefficients.

Let $C$ be a smooth, complete curve over $k$, and let $\delta \subset C$ be a
nonempty finite set of reduced points. Let $J(C)$ be the Jacobian of $C$, and let
$J(C,\delta)$ be the semiabelian variety which represents the functor
\ga{2.1}{ X \mapsto \{(\sL, \phi)\ |\  \sL \text{ line bundle on }
 C \times X, \ {\rm deg} \ \sL|_{C\times k(X)}=0\\
\phi: \sL|_{\delta \times X }\cong \sO_{\delta \times X } \}/\cong . 
\notag}
There is an exact sequence
\eq{2.2}{0 \to T \to J(C,\delta) \to J(C) \to 0
} 
where $T$ is the torus $\Gamma(\delta,\sO^\times)/\Gamma(C,\sO^\times)$. 
By abuse of notation we shall write $J(C, \delta)$ rather than $J(C, \delta)(k)$.
We can
identify the character group $\Hom(T,\G_m)$ with $\Div_\delta^0(C)$, the group of
$0$-cycles of degree $0$ supported on $\delta$. A split subgroup $\Delta \subset
\Div_\delta^0(C)$ corresponds to a quotient $T\surj T_\Delta= T/\ker\Delta$, where
$\ker\Delta\subset T$ is the subtorus killed by all characters in $\Delta$. We may
push out
\eqref{2.2} and define $J(C,\Delta):= J(C,\delta)/\ker\Delta$:
\eq{}{0 \to T_\Delta \to J(C,\Delta) \to J(C) \to 0. 
}
The functor represented by  $J(C,\Delta)$ is the following  quotient of \eqref{2.1}
\ga{2.1a}{ X \mapsto \{(\sL, \phi)\ |\  \sL \text{ line bundle on }X\times C, \ {\rm deg} \ \sL|_{k(X)\times C}=0\\
\phi: \otimes_i
\sL^{\otimes n_i}|_{X\times \{c_i\}} \cong \sO_X \ {\rm for \ all} \ \sum n_i c_i \in \Delta\}.
\notag}
 These trivializations should
be compatible in an evident way with the group law on $\Delta$. 

\begin{lem}\label{lem2.1} We write $H^1(C,\delta)=H^1_c(C\setminus \delta)$. Define
$$H^1(C,\Delta):=(H^1(C,\delta)/\Delta^\perp)\otimes\Q_\ell,$$ 
where
$$\Delta^\perp\otimes\Q_\ell\subset \Q_\ell[\delta]/\Q_\ell \subset
H^1(C,\delta;\Q_\ell)$$  is 
perpendicular to
$\Delta\subset \Div_\delta^0(C)$ under the evident coordinatewise duality. 
The first Chern class $c_1(\sL_\Delta)$ of the
Poincar\'e bundle
$\sL_\Delta$ on $C\times J(C,\Delta)$ lies  in $H^1(C,\Delta) \otimes
H^1(J(C,\Delta))(1)$.
\end{lem}
\begin{proof} Let $I_\delta \subset \sO_{J(C,\Delta)\times C}$ be the ideal of
$J(C,\Delta)\times \delta$. Let $\pi: C \to C'$ be the singular curve obtained
from $C$ by gluing all the points of $\delta$ to a single point $\delta' \in C'$.
Define
$M_\Delta
\subset  (1\times \pi)_*(\sO_{J(C,\Delta)\times C}^\times)/k^\times$ to be the
pullback as indicated: \minCDarrowwidth.2cm
\eq{}{\begin{CD}0 @>>> (1\times \pi)_*(1+I_\delta) @>>> M_\Delta @>>>
(\ker\Delta)_{J(C,\Delta)\times \delta'} @>>> 0
\\ @. @| @VVV @VVV \\
0 @>>> (1\times \pi)_*(1+I_\delta) @>>> \frac{(1\times
\pi)_*\sO_{J(C,\Delta)\times C}^\times}{k^\times} @>>>
\frac{(1\times \pi)_*\sO_{J(C,\Delta)\times\delta}^\times}{k^\times} @>>> 0.
\end{CD}
}
One gets a diagram of Kummer sequences of sheaves on $J(C,\Delta)\times C'$ (Here
$j: C\setminus \delta \inj C$)
\eq{2.5}{\begin{CD}@. 0 @. 0 @. 0 \\
@. @VVV @VVV @VVV \\
0 @>>> (1\times \pi)_*(1\times j)_!\mu_{\ell^n} @>>> (1\times \pi)_*(1+I_\delta)
@>\ell^n >>  (1\times \pi)_*(1+I_\delta) @>>>  0 \\ @. @VVV @VVV @VVV \\
0 @>>>M_{\Delta, \ell^n} @>>> M_\Delta @>\ell^n >> M_\Delta @>>>  0 \\
@. @VVV @VVV @VVV \\
0 @>>> ((\ker\Delta)_{J(C,\Delta)\times \delta'})_{\ell^n} @>>>
(\ker\Delta)_{J(C,\Delta)\times \delta'} @>\ell^n >>
(\ker\Delta)_{J(C,\Delta)\times \delta'} @>>>  0 \\ @. @VVV @VVV @VVV \\ @. 0 @.
0 @. 0
\end{CD}
}
We have $[\sL_\Delta] \in H^1(J(C,\Delta)\times C', M_\Delta)$ and so by the
Kummer coboundary, $c_1(\sL_\Delta) \in \varprojlim_n H^2(J(C,\Delta)\times
C',M_{\Delta, \ell^n})$. But $M_{\Delta, \ell^n} \cong \Z/\ell^n_{J(C,\Delta)}
\boxtimes \psi_{\ell^n}$, where $\psi_{\ell^n}$ fits into an exact sequence of
sheaves on $C'$
\eq{2.6}{0 \to \pi_*j_!\mu_{\ell^n} \to \psi_{\ell^n} \to
(\ker\Delta)_{\delta',\ell^n}
\to 0. }
We can identify $\Delta^\perp\otimes \mu_{\ell^n}$ with $(\ker\Delta)_{\ell^n}$.
the exact cohomology sequence from \eqref{2.6} yields
\eq{}{(\ker\Delta)_{\mu_{\ell^n}}\to H^1(C,\delta;\mu_{\ell^n}) \to H^1(C',
\psi_{\ell^n}) \to 0. 
}
Passing to the limit over $n$, it now follows that we may define $c_1(\sL_\Delta)
\in H^1(J(C,\Delta)) \otimes H^1(C,\Delta)(1)$ as in the statement of the lemma. 
\end{proof}

\begin{lem}\label{lem2.2} Suppose given a morphism $\rho:X \to J(C)$. Let $\Xi$ be
a Cartier divisor on $ C\times X $ representing $\rho$. We assume $\Xi$ is flat 
over
$C$ so we may define a correspondence $\Xi_*: \Div(C) \to \Div(X)$. 
Let $U\subset X$ be  nonemtpy open in $X$. Then there
exists a lifting $\rho_{U, \Delta}: U \to J(C,\Delta)$ of $\rho$ if and only if
$(\Xi|_{ C \times U} )_*(\Delta)
\subset \Div(U)$ consists of principal divisors. The set of such liftings is a
torsor under $\Hom(\Delta, \Gamma(U, \sO_U^\times))$. 
\end{lem}
\begin{proof}Choose a basis $z_i = \sum_j n_{ij}c_j$ for the free abelian group
$\Delta$. Write $\sO_{C\times X}(\Xi)_{z_i\times X} := \otimes_j\sO_{C\times X}(\Xi)^{\otimes n_{ij}}|_{ \{c_j\}\times X}$.
The assumption that $\Xi_*(\Delta)$ consists of principal divisors is precisely
the assumption that all the line bundles 
$\sO_{C\times X}(\Xi)_{z_i\times X}|_U$ are trivial.
The choice of the 
trivializations for a basis of $\Delta$ yields the choice of the desired lifting $\rho_{U,\Delta}$. 
\end{proof}

\begin{lem}Assume $X$ is a smooth variety, and let $\rho : X \to J(C)$ be as
above. Suppose $U \subset X$ is a dense open set such that $\rho|_U$ admits a
lifting $\rho_{U,\Delta} : U \to J(C,\Delta)$. Let $\Div_{X\setminus U}^0(X)$ be the free
abelian group on Cartier divisors supported on $X\setminus U$ which are homologous to $0$
on $X$. Then we get a commutative diagram on cohomology \minCDarrowwidth.2cm
\eq{2.8}{\begin{CD} 0 @>>> H^1(J(C)) @>>> H^1(J(C,\Delta)) @>>> \Delta\otimes
\Q_\ell(-1) @>>> 0 \\
@. @VV\rho^* V @VV\rho^*_{U,\Delta} V @VV a V \\
0 @>>> H^1(X) @>>> H^1(U) @>>> \Div_{X\setminus U}^0(X)\otimes \Q_\ell(-1) @>>> 0.
\end{CD}
}
\end{lem}
\begin{proof} The left hand square is commutative by functoriality. That the
cokernels on the top and bottom row are as indicated follows on the top row from
the Leray spectral sequence for the projection $\pi: J(C,\Delta) \to J(C)$ and on
the bottom from the localization sequence which may be written
\eq{}{0 \to H^1(X) \to H^1(U) \to H^2_{X\setminus U}(X) \to H^2(X).
}
The identification $H^2_{X\setminus U}(X) \cong \Div_{X\setminus U}(X)\otimes \Q_\ell(-1)$ is saying that by purity, the Gysin homomorphism is an isomorphism. 
\end{proof}
\begin{rmk}
(i) Fixing $\rho_{U, \Delta}$ amounts to fixing trivializations of
the restriction  $
\sO_{C\times X}(\Xi)_{z_i\times X}|_U
$ as above. Such trivializations exhibit
$$\sO_{C\times X}(\Xi)_{z_i\times X} \cong \sO_X(D_i)$$ for some divisor $D_i$ with support on $X\setminus U$. The
map labeled $a$ in \eqref{2.8} sends $z_i \mapsto D_i$. \newline\noindent
(ii) The diagram 
\eq{2.10}{\begin{CD} \Delta @>>> J(C) \\
@VV a V @VV \rho^* V \\
\Div_{X\setminus U}^0(X) @>>> \Pic^0(X)
\end{CD}
}
is commutative, where the horizontal arrows are cycle classes. Indeed, both $a$
and $\rho^*$ are defined by the divisor on $C\times X$. Note that $a$ depends on
the choice of $\rho_{U,\Delta}$ but only up to rational equivalence. 
\end{rmk}

Now suppose $X$ is smooth, projective, of dimension $d$. Let $U \subset X$ be
a dense open subset. Write $X\setminus U = D \cup Z$ where $D\subset X$ is a divisor and
${\rm codim}(Z\subset X)\ge 2$. We have $H^1(X\setminus D) \cong H^1(U)$. Since we are
interested in $H^1(U)$, we may assume $U = X\setminus D$ is the
complement of a divisor. 

Let $i: C \inj X$ be a general linear space section of dimension $1$, and let
$\delta = C\cap D$.  We may choose $\rho: X \to J(C)$ such that the composition 
\eq{2.12}{\Pic^0(X) \xrightarrow{i^*} J(C) \xrightarrow{\rho^*} \Pic^0(X)
}
is multiplication by an integer $N \neq 0$. 
Indeed, let $H$ be a very ample line bundle  so that $C$ is the 
$(d-1)$-fold product of general sections of $H$. Intersection with $H$ 
 yields an isogeny $\Pic^0(X)\to {\rm Alb}(X)$, which defines an 
 inverse isogeny  ${\rm Alb}(X)\to \Pic^0(X)$ of degree $N$. We pull back the Poincar\'e bundle from $J(C)\times J(C)$ to $C\times X$ via the 
 composite map $C\times X\to J(C) \times {\rm Alb}(X) \to 
 J(C)\times \Pic^0(X)\to 
 J(C)\times J(C)$, where the first map is the cycle map, the second one is  $1\times $isogeny, the third one is $1\times$ restriction. We define $\sO_{C\times X} (\Xi)$ to be the inverse image of the Poincar\'e bundle.  The morphism  $\rho: X\to J(C)$ is the correspondence $x\mapsto \sO_{C\times X} (\Xi)|_{C\times \{x\}}$
  and does not depend on the choice of the section $\Xi$.

Consider the diagram
\minCDarrowwidth.2cm
\eq{}{\begin{CD} 0 @>>> \Q_\ell[\delta]/\Q_\ell @>>> 
H^1_c(C\setminus \delta) @>>> H^1(C)@>>> 0
\\
@. @VV b V @VVi_* V @VV i_*V \\
0 @>>> \frac{H^{2d-2}(D)}{H^{2d-2}(X)}(d-1) @>>> H^{2d-1}_c(U)(d-1) @>>>
H^{2d-1}(X)(d-1)
\end{CD}
}
Here, the rows are long exact sequences associated to restriction to closed
subsets, and the vertical arrows are Gysin maps. The map $b$ can be described as
follows. The $\Q_\ell$-vector space $H^{2d-2}(D)(d-1)$  
has  basis the
irreducible components of $D$, and $b(x)$ is the basis element $[D_x]$ associated
to the unique component $D_x$ of $D$ containing $x$. We have dual exact sequences
(defining $\Div_D^0(X)$)
\ga{}{0 \to \Div_D^0(X) \to H_{2d-2}(D)(1-d) \to H^2(X)(1) \\
H^{2d-2}(X)(d-1) \to H^{2d-2}(D)(d-1) \to \frac{H^{2d-2}(D)}{H^{2d-2}(X)}(d-1)
\to 0. \notag
}
If we view $\Q_\ell[\delta]$
and $H^{2d-2}(D)(d-1)$ as endowed with symmetric pairings with orthonormal bases
the points $x\in\delta$ and the cohomology classes of irreducible components $D_i
\subset D$, then $b$ is adjoint to the map $D_i \mapsto D_i\cdot\delta$. We
conclude
\begin{lem}\label{lem2.5} Define $\Div_D^0(X)$ to be the $\Q_\ell$-vector space
spanned by divisors on $X$ supported on $D$ and homologous to $0$ on $X$. Define
$\Delta
\subset \Div_\delta^0(C)$ to be the image of $\Div_D^0(X)$ under pullback $i^*$.
Then there is a commutative diagram
\eq{}{\begin{CD} 0 @>>> \Q_\ell[\delta]/\Delta^\perp @>>> H^1_c(C,\Delta) @>>>
H^1(C)@>>> 0
\\
@. @VV b V @VVi_* V @VV i_*V \\
0 @>>> \frac{H^{2d-2}(D)}{H^{2d-2}(X)}(d-1) @>>> H^{2d-1}_c(U)(d-1) @>>>
H^{2d-1}(X)(d-1).
\end{CD}
}
\end{lem}
\begin{proof}The map $b$ is dual to the restriction map $\Div_D^0(X)
\xrightarrow{i^*} \Div_\delta^0(C)$. By definition $\Delta^\perp$ is orthogonal to
the image of $i^*$, i.e. $\Delta^\perp = \ker b$. 
\end{proof}
\begin{lem} \label{lem:lift}
Let $\Delta = i^*(\Div_D^0(X)) \subset \Div_\delta^0(C)$ be as in Lemma
\ref{lem2.5}. Then $\rho$ defined in \eqref{2.12} lifts to some 
$\rho_{U,\Delta}: U\to J(C,\Delta)$.  
\end{lem}
\begin{proof}
The correspondence defined by $\sO_{C\times X}(\Xi)$ in \eqref{2.12}
 carries
$\sO_C(z)$ for $ z\in \Delta = i^*(\Div_D^0(X))$ to line bundles in
$\Pic^0(X)$,  the  classes of which fall in the image of 
$\rho^*i^*(\Div_D^0(X))\equiv N\cdot \Div_D^0(X)$ in  $ \Pic^0(X)$.
 To be more precise, let $D_p$ be a basis for
$\Div_D^0(X)$, and set $z_p = i^*D_p$. This is  a basis of $\Delta$. Then 
$\sO_{C\times X}(\Xi)|_{z_p\times X}=\sO_X(D_p)$. Thus choose $\rho_{U,\Delta}$ in 
 Lemma \ref{lem2.2} using this trivialization on $U$.
\end{proof}
Using the Lemmas \ref{lem2.1}, \ref{lem2.5}, \ref{lem:lift} together with
\eqref{2.8}, we pull back
\ml{}{c_1(\sL_\Delta) \in H^1(C,\Delta) \otimes H^1(J(C,\Delta))(1)
\xrightarrow{i_*\otimes \rho_{U,\Delta}}\\ 
H^{2d-1}_c(U)(d) \otimes H^1(U) \cong  H^1(U)^\vee \otimes H^1(U)
}
and define a correspondence $\Phi: H^1(U) \to H^1(U)$. 
\begin{lem} \label{lem2.7}
The map $\Phi$ is the  multiplication  by  $ N$. 
\end{lem}
\begin{proof}
We consider $\Phi $. It acts on $H^1(U)$, comptibly with the exact sequence
\eq{}{0 \to H^1(X) \to H^1(U) \to \Div_D^0(X)(-1) \to 0
}
By definition of $\rho_{U,\Delta}$, it is equal to $N\cdot {\rm Id}$ on $H^1(X)$ and on 
$\Div_D^0(X)(-1) $. Thus $\Phi -N\cdot {\rm Id}$ is a correspondence from $\Div_D^0(X)(-1)$ to
$H^1(X)$. We use purity in the sense of Deligne. There is no nontrivial correspondence 
$\Div_D^0(X)(-1)\to H^1(X)$. If $k=\C$ and we consider  Betti cohomology, 
$\Div_D^0(X)(-1)$ is pure of weight 2 while $H^1(X)$
is pure of weight 1. If $k$ is the algebraic closure of a finite
 field, we have the same conclusion. Otherwise, all the objects used  are defined over a 
finitely generated field  $k$ over a finite field $k_0 $. By Cebotarev theorem, the Galois group
of $k/k_0$ is generated by Frobenii, so we may make sense of the notion of weight for
$H^1(U)$. We conclude as in the complex case. 
\end{proof}

We now express in terms of cycles the trivialization of $\sO_{C\times X}
(\Xi)_{_p\times X}=\sO_X(D_p)$ used in the proof of Lemma \ref{lem:lift}.  
\begin{thm}With notation as above, there exists a cycle $\Gamma$ on $X\times U$ of
dimension $d= \dim X$ together with rational functions $f_\mu$ on $X$ for each divisor 
$\mu$ homologically equivalent to $0$ on $X$ and supported on $D=X_U$ such that
$pr_{2*}(\Gamma\cdot (\mu\times U)) = (f_\mu)$. The data $(\Gamma, \{f_\mu\})$ define a
class in $H^{2d-1}_c(U)\otimes H^1(U)$ which gives the identity map on $H^1(U)$. 
\end{thm}

We close this section with a comment about K\"unneth projectors 
$\pi^i: H^*(U) \to H^i(U)$ for $i>1$. We consider the somewhat weaker question of the
existence of an algebraic projector when we localize at the generic point of the target,
i.e. we consider $H^*(U) \to H^i(U) \to \varprojlim_{V\subset U} H^i(V)$.  We assume
$U = X\setminus D$ with
$X$ smooth, projective, and
$D$ a Cartier divisor. 
\begin{prop}\label{prop2.9} Let $n< \dim X$ be an integer. Let $Y \subset X$ be a plane
section of dimension $n$ which is general with respect to $D$. Write $\delta = Y\cap D$.
Then the restriction map
\eq{}{H^{n+1}_D(X) \to H^{n+1}_\delta(Y)
}
is injective.
\end{prop}
\begin{proof} Let $d={\rm dim}(X)$. By duality, 
we have to show surjectivity of the Gysin map $H^{n-1}(\delta)\to H^{2d-(n+1)}
(D)(d-n)$.
More generally, one has 
\begin{thm}[P. Deligne]
Let $\sF$ be a $\ell$-adic sheaf 
on $\P^n$, then for $A$ in a 
non-trivial  open subset of the dual projective space 
$(\P^n)^\vee$, the Gysin homomorphism 
\ga{}{H^{i-2}(A, \iota^*\sF)(-1)\to H^i_A(\P^n, \sF),\notag}
where $\iota: A\to \P^n$ is the closed embedding, is an isomorphism
for all  $i$.
In particular, 
if $V \subset \P^n$ is a projective variety, 
then the Gysin homomorphism $H^i(A\cap V)\to H^{i+2}(V)(1)$ is an isomorphism
 for 
$i>{\rm dim}(A\cap V)$ and surjective for $i={\rm dim}(A\cap V)$ for $A$ in a nonempty open subset of $(\P^n)^\vee$. 
\end{thm}
 The proof of the general theorem is written in \cite{Es}, 
Theorem 2.1. Applied to $\sF=a_*\Q_\ell$, where $a: V\to \P^n$ is the projective embedding, it shows that the Gysin  isomorphism $H^i(A\cap V)\to H^{i+2}_{A\cap V}(V)(1)$ is an isomorphism. Then the application follows from Artin's vanishing theorem $H^i(V\setminus (A\cap V))=0$ for $i>{\rm dim}(V)$. 
\end{proof}

Let $L$ be the Lefschetz operator on $H^*(X)$. One of the standard conjectures ($B(X,L)$ in
\cite{Kl2}) is the existence of an algebraic correspondence $\Lambda$ which is a ``weak inverse'' to
$L$. Assume now that this standard conjecture $B$ is true for $X$ and for all smooth
linear space sections $Y \subset X$. The strong Lefschetz theorem implies
that $L^{d-n} : H^n(X)
\xrightarrow{\cong} H^{2d-n}(X)(d-n)$. Assuming $B(X,L)$, $\Lambda^{d-n} = (L^{d-n})^{-1}:H^{2d-n}(X)(d-n)
\cong H^n(X)$. Write $P = \Lambda^{d-n}|_{Y\times X}$. It is easy to check that the composition 

\eq{}{H^n(X) \xrightarrow{i^*} H^n(Y) \xrightarrow{P} H^n(X)
}
is the identity, so $H^n(Y) = \text{Image}(i^*)\oplus \ker(P)$. Consider the diagram
\eq{2.20}{\begin{CD}H^n(X) @>>> H^n(U) @> a >> H^{n+1}_D(X) \\
@VVi^* V @VVV @VV b V \\
H^n(Y) @> d >> H^n(Y\setminus \delta) @>c>> H^{n+1}_\delta(Y)
\end{CD}
}
Define
\eq{}{H^n(Y\setminus \delta)^0 = \{x\in H^n(Y\setminus \delta)\ |\ c(x) \in {\rm Im}(b\circ a)\}.
}
As a consequence of proposition \ref{prop2.9} and \eqref{2.20} we see that $H^n(U) \surj
H^n(Y\setminus \delta)^0/d(\ker(P))$, and the kernel of this map is the image in $H^n(U)$ of elements
$x\in H^n(X)$ such that $i^*x \in \ker(P)\oplus\text{Image}(H^n_\delta(Y) \to H^n(Y))$. For such an
$x$, it will necessarily be the case that $x=P(i^*x)$ is supported on a proper closed
subset of $X$. In particular, for some $V\subset U$ open dense, $P$ will induce a map  
\eq{}{ P_U :H^n(Y\setminus \delta)^0 \to H^n(V). 
}
The map $i^*:H^n(U) \to H^n(Y\setminus \delta)^0$ dualizes to 
$i_*:(H^n(Y\setminus \delta)^0)^\vee  \to
H^{2d-n}_c(U)$, so we may define
\eq{}{(i_*\otimes P_U): H^n(Y\setminus \delta)^0)^\vee \otimes H^n(Y\setminus
\delta)^0 \to H^{2d-n}_c(U)\otimes H^n(V). 
}
Let $T \subset H^n(Y\setminus\delta)^0$ be the subgroup of cohomology classes supported in
codimension $1$. Assuming inductively that we are able to define an algebraic
correspondence on $Y$ which carries a class 
\eq{}{ \gamma \in (H^n(Y\setminus \delta)^0)^\vee \otimes (H^n(Y\setminus
\delta)^0/T) 
}
corresponding to the evident map $H^n(Y\setminus \delta)^0 \surj H^n(Y\setminus
\delta)^0/T$, it would follow since $P_U(T) \subset \ker(H^n(V) \to
\varprojlim_{V\subset U} H^n(V))$ that we could view
\eq{}{(i_*\otimes P_U)(\gamma) \in  H^{2d-n}_c(U)\otimes \varprojlim_{V\subset U} H^n(V).
}
This correspondence would have the desired properties.

One interest in pursuing this line of investigation concerns the Milnor conjecture that 
the
Galois cohomology with $\Z/n\Z$-coefficients prime to the residue characteristic is 
generated as
an algebra by $H^1$. There is a geometric proof of this result in top degree \cite{Bl}, 
so, for example, elements in $H^n(Y\setminus \delta)$ lie in the subalgebra generated by
$H^1$  after localization.
If $P_U$ exists as an algebraic correspondence, then using the existence of a norm in 
Milnor
$K$-theory, one could show that the Milnor conjecture was true for $H^*(k(X), \Z/\ell\Z)$ 
for
almost all $\ell$. (The condition on $\ell$ arises because the standard conjectures only 
make
sense after tensoring with $\Q$.) Here, the idea that cohomology classes in degree $n$ 
might
come by correspondence from an algebraic variety of dimension $\le n$ was suggested to us 
by A.
Beilinson.

\section{Open correspondences}
The aim of this section is to give a simple motiivc proof of the independence of $\ell$ 
or of a complex embedding of a ground field $k$
of the trace for open correspondences. If we assume that $k$ is finite, 
then, as conjectured by Deligne, high Frobenius power twists move the 
correspondence to a general position correspondence and the local factors have been  computed in \cite{De}, \cite{P}, \cite{Z}, \cite{Fu}.  
Surely in this case the simple observations which follow are weaker.

We consider open correspondences. This means the following. Let $X$ be a smooth projective variety of dim $d$ over an algebraically closed field $k$, and let $U\subset X$ be a nontrivial open subvariety,
with complement 
$D=X\setminus U$. One considers  
codim $d$ cycles $\Gamma\subset U\times U$ which have the property that they induce a correspondence $\Gamma_*: H^i(U)\to H^i(U)$ or equivalently $H^i_c(U)\to H^i_c(U)$ for all $i$. Here cohomology is \'etale $\Q_\ell$ cohomology 
or Betti cohomology if $k=\C$ and we denote by $p_i: X\times X\to X$ the two projections. 
We write 
\ga{4.1}{ \Gamma=\sum n_j \Gamma_j}
where $\Gamma_j$ is irreducible, $n_j\in \Z$ and define 
\ga{4.2}{\bar{\Gamma}:=\sum n_j \bar{\Gamma}_j, \ \Gamma_j\subset U\times U}
where $\ \bar{} \ $ is the Zariski closure in $X\times X$. 
\begin{defn} \label{defn4.1}
If $p_2|_{(\Gamma)_j}: \Gamma_j \to U$ is proper for all $j$, or equivalently if 
\ga{4.3}{\bar{\Gamma}_j\cap (D\times X)\subset X\times D \ \forall \ j,}
then one defines
\ga{4.4}{(\Gamma_j)_*: 
H_c^i(U)\xrightarrow{p_2^*} H^i_c(\Gamma_j) \xrightarrow{(p_1)_*}H^i_c(U),} 
and call it the open correspondence defined by $\Gamma_j$. The correspondence
defined by $\Gamma$ is then by definition 
$\Gamma_*=\sum n_j (\Gamma_j)_*$.
 \end{defn}
\begin{rmk}
If we compare this condition to the one yielding to K\"unneth correspondences in sections 2 and 3, then it is much stronger. Indeed, we had $P\subset W\times X$, thus $\bar{\Gamma}=P$ viewed as a cycle in $X\times X$ and 
$P\cap (D\times X)=P\cap (\delta \times X)$ was not empty.  
\end{rmk}
We assume now that we have the assumption as in Definition \ref{defn4.1} and we wish to give conditions under which one can compute the trace of $\Gamma_*$ which is defined by
\ga{4.5}{{\rm Tr}(\Gamma_*):=\sum_{i=0}^{2d}(-1)^i {\rm Tr} 
(\Gamma_*|_{H^i_c(U)}).} 
As it stands, the trace of $\Gamma_*$ depends a priori on $\ell$, or,
 for varieties defined over a field $k$ of characteristic 0, and Betti cohomology taken with respect to a complex embedding $\iota: k\to \C$, it depends on $\iota$. One has
\begin{thm} \label{thm4.3}
Let $X$ be a smooth projective variety of dimension $d$ defined over a field $k$, together with a strict normal crossings divisor $D\subset X$ of open complement $U=X\setminus D$.  Let $\Gamma\subset U\times U$ be  a dimension $d$ cycle defining an open correspondence $\Gamma_*$ on $\ell$-adic cohomology or Betti cohomology as in Definition \ref{defn4.1}. Then ${\rm Tr}(\Gamma_*)$ does not depend on $\ell$ in $\ell$-adic cohomology or on the complex embedding of $k$ in Betti cohomology. 
\end{thm}
\begin{proof}
We use the
relative motivic cohomology
$H^{2d}_M(X \times U, D\times U,\Z(d))$,  as defined in
\cite{Le}, chapter 4, 2.2 and p. 209. 
The  group $H^m_M(X \times U, D\times U,\Z(n))$ is the homology
$H_{2n-m}(\sZ^n(X\times U, D\times U,*))$, where
 $\sZ^n(X \times U, D\times U,*)$ is the single
complex associated to the double
 higher Chow cycle complex
\eq{4.6}{\begin{CD}
\cdots &  &\cdots & &\cdots\\
@V\partial VV @V\partial VV @VVV\\
\sZ^n(X\times U, 1)@>{\rm rest}>> \sZ^n(D^{(1)} \times U, 1)@>{\rm rest}>>
\sZ^n(D^{(2)}\times U,1)\\
@V\partial VV @V\partial VV @VVV\\
\sZ^n(X\times U, 0)@>{\rm rest}>> \sZ^n(D^{(1)} \times U, 0)@>{\rm rest}>>
\sZ^n(D^{(2)}\times U,0).
\end{CD}}
Here $D^{(a)}$ is the normalization of all the strata of codimension $a$,
 $\sZ^n(D^{(a)}\times U, b)$ is a group of cycles on
$D^{(a)}\times U\times S^b$ where
$S^\bullet$ is the cosimplicial scheme $S^n = \Spec(k[t_0,\dotsc,t_n]/
(\sum_{i=0}^n t_i -1))$ with
face maps $S^n \inj S^{n+1}$ defined by $t_i=0$.
More precisely, $\sZ^n(D^{(a)}\times U,
b)$ is generated by the codimension $n$ subvarieties
$Z\subset D^{(a)}\times U \times S^b$ such
that, for each face $F$ of
$S^b$, and each irreducible component $F'\subset D^{(a)}$
of the strata of $D$ we have
$
\codim_{F'\times U\times F}(Z\cap (F'\times U\times F))\ge n.
$
The horizontal restriction maps
are the intersection with the smaller strata, the vertical
$\partial$'s are the boundary maps.

This relative motivic cohomology acts as correspondences on
$H^*_c(U)$, where $H^*_c(U)$  is $\ell$-adic or Betti cohomology. 

Let us write $\Gamma=\sum n_j \Gamma_j$. By \eqref{4.3}, one has $\Gamma_j\subset X\times U$ 
closed with $\Gamma_j\cap (D\times U)=\emptyset$, thus in particular, $\Gamma \in \sZ^d(X\times
U,0)$ with ${\rm rest}(\Gamma)=0$ in $\sZ^d(D^{(1)} \times U, 0)$, thus it defines a class
\ga{4.7}{[\Gamma]\in  H^{2d}(X\times U, D\times U, d).}
Similarly, we consider the restriction $\Delta_U\subset U\times X$ of 
the diagonal $\Delta\subset X\times X$. This defines a class in $\sZ^d(U\times X,0)$. As  ${\rm
rest}(\Gamma_U)=0$ in $\sZ^d(U\times D^{(1)}, 0)$,  it defines a class
\ga{4.8}{[\Delta_U]\in H^{2d}(U\times X, U\times D, d).}
We want to pair $[\Gamma]$ with $[\Delta_U]$. We argue using M. Levine's work.
Let $Y$ be a $N$-dimensional 
smooth projective defined over $k$, with two strict 
normal crossings divisors $A,B$ so that  $A+B$ is  a strict normal corssings divisor.  
By \cite{Le}, Chapter IV, lemma 2.3.5 and lemma 2.3.6, 
the motive $M(Y\setminus A, B \setminus B\cap A)$ is dual to the motive 
$M(Y\setminus B, A \setminus B\cap A)$. It yields a cup product 
\ga{4.9}{H^a_M(Y\setminus A, B\setminus B\cap A,b)\times 
H^{2N-a}_M(Y\setminus B, A\setminus A\cap B, N-b) \to \Z.}  
This cup product is compatible with the cup product in $\ell$-adic or Betti cohomology. 
We apply this to
\ga{4.10}{Y=X\times X, \ A=D\times X, \ B=X\times D}
so we can cup $[\Delta_U]\in H^{2d}_M(Y\setminus  A, B\setminus B\cap A, d)$ and  
$[\Gamma]\in
H^{2d}_M(Y\setminus B, A\setminus A\cap B, d)$
\ga{4.11}{[\Delta_U]\cup [\Gamma]\in \Z.}
The theorem is then the consequence of the following proposition.
\begin{prop} \label{prop4.4}
$${\rm Tr}(\Gamma_*)=[\Delta_U]\cup [\Gamma]\in \Z.$$
\end{prop} 
\begin{proof}
By the compatibility
of the cup product \eqref{4.9} with cohomology, we just have to prove the proposition with 
$[\Delta_U]$ and $[\Gamma]$ replaced by their classes 
${\rm cl}(\Delta_U)$ and ${\rm cl}(\Gamma)$ in cohomology.  We may assume that  $\Gamma\subset X
\times U$ is irreducible. On the other hand, the map
\ga{4.12}{ H^i_c(\Gamma)\xrightarrow{p_{1*}}
H^i_c(U)}
factors through
\ga{4.13}{H^i_c(\Gamma) \xrightarrow{{\rm Gysin}} H^{i+2d}_{\Gamma, c}(X\times U,d)=
H^{i+2d}_{\Gamma,c}(U\times U,d) \\
\to H^{i+2d}_{c}(U\times U,d) \to H^{i+2d}_c(U\times X, d)
\xrightarrow{p_{1*}} H^i_c(U)\notag
}
so the correspondence
\ga{4.14}{\Gamma_*: H^i_c(U)\xrightarrow{p_2^*}H^i_c(\Gamma)\xrightarrow{p_{1*}} H^i_c(U)} is just
\ga{4.15}{H^i_c(U)\cup ({\rm cl}(\Gamma)\in H^{2d}(X\times U, D\times U, d))
\subset 
H^{i+2d}_c(U \times X, d)}
followed by $p_{1*}$. 
Now we argue as in the classical case. Let $e^i_a$ be a basis of $H^i_c(U)$, and $(e^i_a)^\vee$ 
be its dual basis in $H^{2d-i}(U)(d)$.  Write ${\rm cl}(\Gamma)=\sum_i \sum_a f^i_a \otimes
(e^i_a)^\vee, f^i_a=
\sum f^i_{ab}e^i_b\in H^i_c(U)$. So $\Gamma_*(e^i_a)=\sum_b f^i_{ab}$,
 and 
${\rm Tr}(\Gamma_*)=\sum (-1)^i \sum_a \sum_b f^i_{ab}$. On the other hand, one has  ${\rm
cl}(\Delta_U)=\sum_i\sum_a (e^i_a)^\vee \otimes (e^i_a)$. Thus 
 ${\rm cl}(\Delta_U\cup \Gamma)=\sum (-1)^i \sum_a \sum_b f^i_{ab}$.

\end{proof}
The proposition finishes the proof of  the theorem.
\end{proof}
The rest of the section is devoted to giving a concrete expression for \eqref{4.5}  under
stronger geometric assumptions on $\Gamma$.

The condition \eqref{4.3} allows to define the embedding
\ga{4.16}{
\Gamma_j=(p_2|_{\Gamma_j})^{-1}(U)=\bar{\Gamma}_j\setminus 
\bar{\Gamma}_j\cdot (X\times D)\\
\xrightarrow{\iota \ \subset} (\Gamma_j)':=(p_1|_{\Gamma_j})^{-1}(U)=
\bar{\Gamma}_j\setminus 
\bar{\Gamma}_j\cdot (D\times X).
\notag}
By definition
\ga{4.17}{p_1|_{(\Gamma_j)'} \ {\rm is \ proper}.}
\begin{lem} \label{lem4.5}
For all $j$, one has a factorization
\ga{}{(\Gamma_j)_*: H^i_c(U)\xrightarrow{p_2^*} H^i_c(\Gamma_j)\xrightarrow{\iota_*}
H^i_c((\Gamma_j)') \xrightarrow{(p_1)_*} H^i_c(U).\notag} 

\end{lem}
\begin{proof}
By definition one has a factorization
\ga{4.18}{
\begin{CD}
\Gamma_j @> \iota >> (\Gamma_j)' @> \mu >> \bar{\Gamma}_j\\
@V p_1|_{\Gamma_j} VV @ VV p_1|_{(\Gamma_j)'}V @VV p_1|_{\bar{\Gamma}_j} V \\
U @ >>=> U  @>>j > X
\end{CD}
}
Setting $\lambda= \mu \circ \iota$, the lemma then just says that on 
$p_1$-acyclic locally constant sheaves, one has a factorization
$ (p_1|_{\bar{\Gamma}_j})_!\lambda_!\to  (p_1|_{\bar{\Gamma}_j})_!\mu_!
\xrightarrow{{\rm Tr}} j_!$.

\end{proof}
\begin{cor} \label{cor4.6}
One has a commutative diagram
\ga{}{\begin{CD}
H^i_c(U) @>>> H^i(X) \\
@V p_2^* VV  @V p_2^* VV   \\
H^i_c(\Gamma_j) @>>> H^i(\bar{\Gamma}_j) \\
@V \iota_* VV @V = VV 
\\
H^i_c((\Gamma_j)') @>>> H^i(\bar{\Gamma}_j) 
\\
@V p_{1*} VV @V p_{1*} VV  \\
H^i_c(U) @>>> H^i(X) 
\end{CD}\notag}
 \end{cor}
\begin{proof}
With the notations of \eqref{4.18}, the commutativity of the two upper squares simply means on 
constant sheaves $p_2j_!=\lambda_!p_2^*
\to \mu_!p_2^*=\mu_!$. Thus the assertion is on the lower square. 
By \eqref{4.17}, one has $R((p_1)_*)|_{\Gamma'_j}=R((p_1)_!)|_{\Gamma'_j}$
thus $R((p_1)_*)|_{\bar{\Gamma}_j}\mu_!= R((p_1)_!)|_{\bar{\Gamma}_j}\mu_!=
j_!=R((p_1)_!)|_{\Gamma'_j}=j_! R((p_1)_*)|_{\Gamma'_j}$. This shows the commutativity. 

\end{proof}
We remark that in Corollary \ref{cor4.6}, 
the Grothendieck-Lefschetz trace formula allows to compute the trace of the correspondence 
$\bar{\Gamma}_*$ on $X$
\ga{4.19}{
{\rm Tr}(\bar{\Gamma}_*)={\rm deg} (\bar{\Gamma}\cdot \Delta_X)
}
Thus, in the corollary, 
 we would like to complete the commutative diagram in an exact sequence of commutative diagrams, 
so that we can apply  the trace formula on all the terms but the one we seek. 
In the sequel, we give a strong geometric condition under
 which it is possible.

\begin{defn} \label{defn4.7}
We assume that $D$ is a strict normal crossings divisor. 
The dim $d$ cycle $\Gamma \subset U\times U$
is said to be in good position with respect to $D\times X$ 
if the following two conditions are fulfilled.
\begin{itemize}
\item[i)] Each
$\bar{\Gamma}_j$ cuts each stratum $D_I\times X$
in codim $\ge d$,  where 
$\ D_I=D_{i_1}\cap \ldots \cap D_{i_r}$ for $ \  I=\{i_1,\ldots, i_r\}$  with $|I|=r$. 
\item[ii)]
\ga{}{
\bar{\Gamma}_j \cap (D_I\times X)\subset D_I\times D_I\notag
}
set theoretically. 
\end{itemize}
In this case, for all $I$ we  define the cycles 
\ga{4.20}{
Z_{jI}=\bar{\Gamma}_j\cdot (D_I\times X)\subset D_I\times D_I
}
\end{defn}
Let us be more precise. We denote motivic cohomology by $\sH^a(b)$. We drop the subscript $_j$, 
thus $\Gamma=\Gamma_j$. This defines 
\ga{4.21}{Z_I=\sum m_{I,a} Z_{I,a}} 
as in \eqref{4.10}, where the $Z_{I,a}$ are the reduced irreducible components of $Z_I$.
One has the Gysin isomorphisms
\ga{4.22}{\oplus_a \Q_\ell[Z_{I,a}] \xrightarrow{\cong} 
H^{2(d-r)}_{Z_I}(D_I\times D_I, d-r)\xrightarrow{\cong}
H^{2d}_{Z_{I}}(D_I\times X, d)
}
This yields the commutative diagram
 \ga{4.23}{
\begin{CD}
\oplus_a\Q\cdot [Z_{I,a}] @>\otimes_{\Q} \Q_\ell >>  
\oplus_a \Q_\ell \cdot [Z_{I,a}]\\
@V\cong VV @VV\cong V \\
\sH^{2(d-r)}_{Z_I}(D_I\times D_I, d-r ) @>\otimes_{\Q} \Q_\ell>>  
H^{2(d-r)}_{Z_I}(D_I\times D_I, d-r )  \\
@V\cong VV  @VV\cong V \\ 
\sH^{2d}_{Z_I}
(D_I\times X, d)@>\otimes_{\Q} \Q_\ell >> H^{2d}_{Z_I}(D_I\times X, d)
 \end{CD}
}
So we conclude that $Z_I$ is a well defined cycle  
\ga{4.24}{
Z_I\in \sH^{2(d-r)}_{Z_I}(D_I\times D_I, d-r ) \to \sH^{2(d-r)}(D_I\times D_I, d-r ),
}
 which defines a correspondence
\ga{4.25}{
(Z_I)_*=\sum_a m_{I,a}(Z_{I,a})_*: H^i(D_I)\to H^i(D_I),}
by the formula
\ga{4.26}{(Z_{I,a})_*: H^i(D_I)\xrightarrow{p_2^*} H^i(Z_{I,a}) \xrightarrow{(p_1)_*} H^i(D_I).}
It can be described as follows. Let $(Z_{I,a})_i\in H^{2d-2r-i}(D_I)\otimes H^i(D_I)(d-r)$ be the 
K\"unneth component of $Z_{I,a}$. Then one has
\ga{4.27}{
(Z_{I,a})_*: H^i(D_I)\xrightarrow{p_2^*} H^i(D_I)\otimes H^0(D_I)
\xrightarrow{\cup (Z_{I,a})_i}\\ H^i(D_I)\cup  H^{2d-2r-i}(D_I)(d-r)\otimes H^i(D_I)
 \xrightarrow{{\rm Tr}\otimes 1} H^i(D_I). \notag
}
We denote by $\Delta_I\subset D_I\times D_I$ the diagonal. By the 
Grothendieck-Lefschetz trace formula, one has
\ga{4.28}{{\rm Tr} ((Z_I)_*)={\rm deg}( 
Z_I \cdot \Delta_I).}
\begin{lem} \label{lem4.8}
Let $J\subset I$ with $|J|=s, |I|=r$. One has a commutative diagram
\ga{}{\begin{CD}
H^i(D_I)@>>{\rm rest}> H^i(D_J)\\
@V (Z_I)_* VV @VV (Z_J)_* V \\
H^i(D_J) @>>{\rm rest}> H^i(D_J)
\end{CD} \notag}
\end{lem}
\begin{proof}
We have $Z_I=\bar{\Gamma}\cdot (D_I\times X), 
Z_J=\bar{\Gamma}\cdot (D_J\times X)=Z_I\cdot (D_J\times X)$.  
One has the diagram 
\ga{4.29}{\begin{CD}
Z_J@>>\subset > D_J\times D_J \\
@. @VV\subset V\\ 
@. D_J\times D_I @>> \subset > D_J\times D_I\\
@. @. @AA \subset A\\
@. @. Z_I
\end{CD}
}
yielding the commutative diagram 
\ga{4.30}{\begin{CD}
H^{2d-2r-i}(D_I)\otimes H^i(D_I)(d-r)@>>1\otimes {\rm rest}> 
 H^{2d-2r-i}(D_I)\otimes H^i(D_J)(d-r)\\
@. @AA{\rm Gysin}\otimes 1 A\\
@. H^{2d-2s-i}(D_J)\otimes H^i(D_J)(d-s)
\end{CD}
}
Via this diagram, one has
\ga{4.31}{ (1\otimes {\rm rest})((Z_I)_i)= ({\rm Gysin}\otimes 1)((Z_J)_i)}
Here the notation $Z_I$ means the sum aover all components with multiplicities, and similarly for 
the K\"unneth components.  Thus we  apply the description \eqref{4.27} to conclude. 
\end{proof}
In order to have a unified notation, we set for $|I|=0$
\ga{4.32}{Z_I=\bar{\Gamma}, \ \Delta_I=\Delta_X .}
One has
\begin{thm} \label{thm4.9}
Let $X$ be a smooth projective variety over an algebraically closed field $k$, $D\subset X$  be a
strict normal crossing divisor, with complement $U=X\setminus D$, and $\Gamma \subset U\times U$
be a dim $d$ correspondence, with $p_2|_\Gamma: \Gamma \to U$ proper, and  in good position with
respect to $D\times X$ in the sense of Definition \ref{defn4.7}. Then one has
\ga{}{
{\rm Tr}(\Gamma_*)=\sum_{r=0}^d (-1)^r \sum_{|I|=r} 
{\rm deg}( 
Z_I \cdot \Delta_I)
.\notag
} 
\end{thm}  
\begin{proof}
One has the Mayer-Vietoris spectral sequence
\ga{4.33}{
E_1^{a,b}=\oplus_{|I|=a}H^b(D_I) \Rightarrow H^{a+b}_c(U)
}
On the other hand, the correspondence $\Gamma_*$ on $H^\bullet_c(U)$ is compatible with the 
correspondence $\bar{\Gamma}_*$ on $H^\bullet(X)$ by  Corollary \ref{cor4.6}, and with the
correspondence $(Z_I)_*$ in $H^\bullet(D_I)$ by Lemma 
\ref{lem4.8}.  Thus $\Gamma_*$ acts on the whole spectral sequence via $\bar{\Gamma}_*, (Z_I)_*$. 
One obtains
\ga{4.34}{
{\rm Tr}(\Gamma_*)=\sum_{r=0}^d (-1)^{r} \sum_{|I|=r} 
{\rm Tr}((Z_I)_*).
}
We conclude with \eqref{4.19}, \eqref{4.28}.
\end{proof}
We now give a scheme-theoretic condition under which the expression given in 
Theorem \ref{thm4.9} depends only on local contributions in $U$. This condition is inspired by 
\cite{KS}, Lemma 2.3.1.  
\begin{defn} \label{defn4.10} 
We assume that $D$ is a strict normal crossings divisor. 
The dim $d$ cycle $\Gamma \subset U\times U$
is said to be in scheme theoretic good position with respect to $D\times X$ 
if it is in good position in the sense of Definition \ref{defn4.7} and ii) is replaced by 
\begin{itemize}
\item[ii')]
\ga{}{
\bar{\Gamma}_j \cap (D_I\times X)\subset D_I\times D_I\notag
}
scheme theoretically. 
\end{itemize}

\end{defn}
\begin{prop} \label{prop4.11}
Let $X$ be a smooth projective variety over an algebraically closed field $k$, $D\subset X$  be a
strict normal crossing divisor, with complement $U=X\setminus D$, and $\Gamma \subset U\times U$
be a dim $d$ correspondence, with $p_2|_\Gamma: \Gamma \to U$ proper, and  in scheme theoretic 
good position with respect to $D\times X$ in the sense of Definition \ref{defn4.10}. We assume
moreover that $
\bar{\Gamma}_j$ and $\Delta|_U$ cut transversally.  Then ${\rm Tr}(\Gamma_*)={\rm deg}
(\Delta_U\cdot \Gamma)$.  
\end{prop}
\begin{proof}
Due to the good position assumption, all intersection multiplicities are $1$ and the contributions
lying on 
$(D\times X) \cup  (X\times D)$ cancel in Theorem \ref{thm4.9}.

\end{proof}
\begin{ex}One case where the conditions of proposition \ref{prop4.11} hold is in
characteristic $p$ when $\Gamma$ is the graph of frobenius. In this case, of course, the result
is known by other methods.
\end{ex}
\begin{ex}
This example is inspired by \cite{KS}, Remark 2.3.6. We take $X=\P^1, D=\{\infty\}, U=\A^1, 
\Gamma=\Gamma_{pq}=\{x^p-y^q=0\}\subset \A^1\times \A^1$.
 Then $\Gamma$ defines an open correspondence with ${\rm Tr}(\Gamma_*)=p$.  
On the other hand, one has
\ga{4.35}{
{\rm deg}(\Gamma_{pq}\cdot \Delta_U)={\rm dim} k[t]/(t^p-t^q)=\\
\begin{cases} 
{\rm max}(p,q) & {\rm if} \ p\neq q\\
\infty & {\rm if} \ p=q
\end{cases}\notag}
and one has
\ga{4.36}{{\rm deg} (\bar{\Gamma}_{pq}\cdot (\infty\times \P^1))=q, \
{\rm deg} (\bar{\Gamma}_{pq}\cdot (\P^1 \times \infty))=p.}
Thus $\Gamma_{pq}$ is in scheme theoretic good position with respect to $\infty \times
\P^1$  if and only if $p>q$, and is always  in good position with respect to $\infty
\times \P^1$. Since
 ${\rm deg}(\Delta_{\P^1}\cdot \bar{\Gamma}_{pq})={\rm deg}(\sO(1,1)\cdot \sO(p,q))=p+q$, we  see
exactly how the formula of Theorem 
\ref{thm4.9} works both in Theorem \ref{thm4.9} and in Proposition \ref{prop4.11}. 
\end{ex}

\newpage
\bibliographystyle{plain}
\renewcommand\refname{References}

\end{document}